\UseRawInputEncoding
\documentclass[oneside, 11pt]{amsart}
\usepackage{amsmath}
\usepackage{amssymb}
\usepackage{dsfont}
\usepackage{color}
\usepackage[usenames,dvipsnames,x11names,svgnames]{xcolor}
\usepackage[colorlinks=true,linkcolor=NavyBlue,citecolor=DarkGreen, urlcolor=Blue]{hyperref}
\usepackage{amsthm}
\usepackage{amscd}
\usepackage{geometry}
\usepackage{mathrsfs}
\usepackage{mathtools}

\newtheorem{theorem}{Theorem}

\newtheorem{lemma}{Lemma}
\newtheorem{corollary}{Corollary}

\theoremstyle{definition}
\newtheorem{definition}{Definition}
\newtheorem*{remark}{Remark}

\newcommand{\be}{\begin{equation*}}
\newcommand{\ee}{\end{equation*} }

\newcommand{\ben}{\begin{equation}}
\newcommand{\een}{\end{equation} }

\newcommand{\bs}{\begin{split}}
\newcommand{\es}{\end{split}}

\newcommand{\bmu}{\begin{multline*}}
\newcommand{\emu}{\end{multline*}}

\newcommand{\bmun}{\begin{multline}}
\newcommand{\emun}{\end{multline}}

\begin{document}



\title[]{A new class of the entire function of order one: a case study}

\author{Xiao-Jun Yang$^{1,2}$}

\email{dyangxiaojun@163.com; xjyang@cumt.edu.cn}

\address{$^{1}$ School of Mathematics,and State Key Laboratory for Geo-Mechanics and
Deep Underground Engineering, China University of Mining and Technology, Xuzhou 221116, China}
\address{$^{2}$ Department of Mathematics, Faculty of Science, King Abdulaziz University P.O. Box 80257, Jeddah 21589, Saudi Arabia}

\begin{abstract}
In this article, a new class of the entire function of order one, expressed
by the series and product representations with the real positive
coefficients and complex zeros, is investigated for the first time. The
entire function on the critical line deduces an even entire function of
order one. It is proved that the real part of the complex zeros is equal to
the critical line. An equivalent representation theorem is obtained to set
up the sufficient conditions for the critical line for the entire function.
As a typical example, the critical line for the special hyperbolic cosine
function obtained by the present theorem agrees with the result of Euler.
We also discover the new products of the hyperbolic cosine
and sinc functions.
\end{abstract}

\maketitle
\tableofcontents

\section{Introduction} {\label{sec:1}}

The distribution of the zeros of the entire functions is one of important
topics in the study of the analytic functions of complex variable \cite{1,2}.
There exists a class of the entire function with their series and product
representations. This idea for this class of the entire functions was proposed
by Laguerre \cite{3} and P\'{o}lya \cite{4}. Let $\mathcal{C} $, $\mathcal{R}$ and $\mathcal{A}$ denote
the sets of the complex, real and integer numbers and let $\mathbf{\Xi} =\mathcal{C}
\backslash \left\{ 0 \right\}$, $\mathbf{H}=\mathcal{R} \backslash \left\{ 0 \right\}$
and $\mathbf{\Pi} =\mathcal{A}\cup \left\{ 0 \right\}$. Suppose that $m\in \mathbf{\Pi} $ and
$k\in \mathcal{A}$. This is so-called Laguerre-P\'{o}lya ($\mathcal{L-P}$) class of the real
entire function $\mathfrak{G}\left( t \right)$, which have the form (see, for instance, \cite{5,6,7,8} and the references therein)
\begin{equation}
\label{eq1}
\mathfrak{G}\left( t \right)=\sum\limits_{m=0}^\infty {\frac{\gamma _m }{m!}t^m} =\alpha
t^n\exp \left( {-bt^2+at} \right)\prod\limits_{k=1}^\varpi {\left(
{1+\frac{t}{t_k }} \right)\exp \left( {-\frac{t}{t_k }} \right)} ,
\end{equation}
where $\alpha ,a,t_k \in \mathcal{R}$, $t_k \ne 0$, $b\ge 0$, $n\in \mathbf{\Pi} $, $0\le
\varpi \le \infty $, and $\sum\limits_{k=1}^\infty {t_k^{-2} } $ is convergent.

A family of the entire functions in the $\mathcal{L-P}$ class has an
important relation with the well-known conjectures \cite{8,9}, for instance, the
Riemann hypothesis \cite{10}, the conjecture of P\'{o}lya \cite{11,12}, the de
Bruijn-Newman constant \cite{13,14}, and the conjecture of Karlin \cite{15}. The zeros
of the successive derivatives of the even $\mathcal{L-P}$ functions were
considered by Shen in \cite{16}. A generalized version of the $\mathcal{L-P}$
class of the entire functions was proposed by Su\'{a}rez in \cite{17}. The
Laguerre polynomials were developed by Dimitrov and Cheikh as the Jensen
polynomials in sense of the $\mathcal{L-P}$ entire functions \cite{18}. The
partial theta function was investigated by Bohdanov and Vishnyakova in the
$\mathcal{L-P}$ class \cite{19}.

In the working paper, we now consider a new class of the entire function as
follows:
\begin{definition}
Let $\ell \in \mathbf{H}$ and $s\in \mathcal{C}$. A entire function $\mathcal{H} \left( s
\right)$ of order $\nu =1$, given by the series
\begin{equation}
\label{eq2}
\mathcal{H} \left( s \right)=\sum\limits_{m=0}^\infty {\Omega _m \left( {s-\ell }
\right)^{2m}} ,
\end{equation}
is said to be in the class $\mathcal{J}$, written $\mathcal{H} \in \mathcal{J}$, if $\mathcal{H}
\left( s \right)$ can be expressed in the product
\begin{equation}
\label{eq3}
\mathcal{H} \left( s \right)=\mathcal{H} \left( 0 \right)\prod\limits_{\rho _k } {\left(
{1-\frac{s}{\rho _k }} \right)} ,
\end{equation}
where $\rho _k \in \mathbf{\Xi}$ take over the zeros of $\mathcal{H} \left( s \right)$,
$\Omega _m >0$ are the real coefficients for $\mathcal{H} \left( s \right)$,
\begin{equation}
\label{eq4}
\hbar =\sum\limits_{k=1}^\infty {\left| {\rho _k } \right|^{-2}}
\end{equation}
is convergent and
\begin{equation}
\label{eq5}
\mathcal{H} \left( 0 \right) \ne 0.
\end{equation}
\end{definition}
Here, we allow that $\alpha _m $ can be replaced by the Taylor series of the
entire function $\mathcal{H} \left( s \right)$. It is clear that
$\mathcal{H} \in \mathcal{J}$ has the infinity of the complex zeros,
$\mathcal{H} \left( 0 \right)=\sum\limits_{m=0}^\infty {\Omega _m \ell^{2m}}>0$ and
$\mathcal{H} \left( \ell \right)=\Omega _0>0$. In fact, two special cases
of the class $\mathcal{J}$ have been considered in \cite{20,21}. The behaviors for
the complex zeros of the class of the entire function
considered in the paper have not been investigated. By the motivation of the
idea, the main aim of the present paper is to prove the following:

\begin{theorem}
\label{TH1}
If $\mathcal{H} \in \mathcal{J}$, then all zeros of $\mathcal{H} \left( s \right)$ lie on
the critical line $Re\left( s \right)=\ell $.
\end{theorem}
Let $Re(s)$ and $Im(s)$ be the real and imaginary parts of a complex variable $s \in \mathcal{C}$.
To prove Theorem \ref{TH1}, we establish a class of the function $\mathcal{H}
\left( s \right)$, given as
\begin{equation}
\label{eq6}
\mathcal{H} \left( s \right)=\sum\limits_{m=0}^\infty {\Omega _m \left( {s-\ell }
\right)^{2m}} =\mathcal{H} \left( \ell \right)\prod\limits_{Im\left( {\rho _k }
\right)>0} {\left[ {1-\left( {\frac{s-\ell }{\rho _k -\ell }} \right)^2}
\right]} .
\end{equation}
With the aid of (\ref{eq6}), we may derive that
\begin{equation}
\label{eq7}
\psi \left( \vartheta \right)=\sum\limits_{m=0}^\infty {\Omega _m \vartheta
^{2m}} =\mathcal{H} \left( \ell \right)\prod\limits_{Im\left( {\rho _k } \right)>0}
{\left[ {1-\frac{\vartheta ^2}{\left( {\rho _k -\ell } \right)^2}} \right]},
\end{equation}
where $\vartheta \in \mathcal{C}$.
Let $\overline {\rho _k } $ be the complex conjugate of the complex zeros
$\rho _k $ and let $\overline \vartheta \in \mathcal{C}$ be the complex conjugate
of a complex variable $\vartheta \in \mathcal{C}$.
Adopting (\ref{eq7}) to find the complex conjugate $\overline {\psi \left( \vartheta
\right)} $ of the function $\psi \left( \vartheta \right)$, we present
\begin{equation}
\label{eq8}
\psi \left( {\overline \vartheta } \right)=\sum\limits_{m=0}^\infty {\Omega _m
 \overline \vartheta ^{2m}} =\mathcal{H} \left( \ell
\right)\prod\limits_{Im\left( {\rho _k } \right)>0} {\left[
{1-\frac{\overline \vartheta ^2}{\left( {\rho _k -\ell } \right)^2}}
\right]} =\mathcal{H} \left( \ell \right)\prod\limits_{Im\left( {\rho _k }
\right)>0} {\left[ {1-\frac{\overline \vartheta ^2}{\left( {\overline {\rho
_k } -\ell } \right)^2}} \right]} .
\end{equation}
to obtain the convergent series
\begin{equation}
\label{eq9}
\sum\limits_{k=1}^\infty {\frac{1}{\left( {\rho _k -\ell } \right)^2}}
=\sum\limits_{k=1}^\infty {\frac{1}{\left( {\overline {\rho _k } -\ell }
\right)^2}}.
\end{equation}

Making use of (\ref{eq9}), we obtain
\begin{equation}
\label{eq10}
\left( {\rho _k -\ell } \right)^2=\left( {\overline {\rho _k } -\ell }
\right)^2
\end{equation}
to reduce to Theorem \ref{TH1}.

The structure of the paper is designed as follows. In Section 2 we address
the functional equation, products and symmetric lines for the entire
function $\mathcal{H} \left( s \right)$. In Section 2 we give the detailed proof of
Theorem \ref{TH1}. In Section 3 we propose an equivalent representation of Theorem
\ref{TH1}. In Section 4 we investigate the critical line of the
special hyperbolic cosine function $F\left( s
\right)=cosh\left( {s-6} \right)$. Finally, we draw the conclusion in
Section 5.

\section{New results} {\label{sec:2}}

\subsection{The functional equation for $\mathcal{H}
\left( s \right)$}
To be begin with we present the functional equation of $\mathcal{H}
\left( s \right)$.

\begin{lemma}
\label{LE1}
Assume that $s\in \mathcal{C}$ and $\mathcal{H} \in \mathcal{J}$. Then we have
\begin{equation}
\label{eq11}
\mathcal{H} \left( s \right)=\mathcal{H} \left( {2\ell -s} \right).
\end{equation}
\end{lemma}

\begin{proof}
By using (\ref{eq2}), we consider
\begin{equation}
\label{eq12}
 \begin{aligned}
\mathcal{H} \left( {2\ell -s} \right)
=\sum\limits_{m=0}^\infty {\Omega _m \left[ {\left( {2\ell -s} \right)-\ell }
\right]^{2m}}
 =\sum\limits_{m=0}^\infty {\Omega _m \left( {\ell -s} \right)^{2m}}
=\sum\limits_{m=0}^\infty {\Omega _m \left( {s-\ell } \right)^{2m}} ,
 \end{aligned}
\end{equation}
which reduce to (\ref{eq11}).

We thus complete the proof of Lemma \ref{LE1}.
\end{proof}

\subsection{The different products for $\mathcal{H} \left( s \right)$}

\begin{lemma}
\label{LE2}
Let $s\in \mathcal{C}$ and $\mathcal{H} \in \mathcal{J}$. Then there exist any constant $\beta
\in \mathcal{C}$ with $\beta \ne \rho _k $ such that
\begin{equation}
\label{eq13}
\mathcal{H} \left( s \right)=\mathcal{H} \left( \beta \right)\prod\limits_{\rho _k }
{\left( {1-\frac{s-\beta }{\rho _k -\beta }} \right)}.
\end{equation}
\end{lemma}

\begin{proof}
Since $\mathcal{H} \in \mathcal{J}$, we have
\begin{equation}
\label{eq14}
 \begin{aligned}
\mathcal{H} \left( s \right)&=\mathcal{H} \left( 0 \right)\prod\limits_{\rho _k } {\left(
{1-\frac{s}{\rho _k }} \right)} =\mathcal{H} \left( 0 \right)\prod\limits_{\rho _k
} {\frac{\rho _k -s}{\rho _k }} \\
&=\mathcal{H} \left( 0 \right)\prod\limits_{\rho _k } {\left( {\frac{\rho _k -\beta
}{\rho _k -\beta }\cdot \frac{\rho _k -s}{\rho _k }} \right)} \\
&=\mathcal{H} \left( 0
\right)\prod\limits_{\rho _k } {\left( {\frac{\rho _k -\beta }{\rho _k
}\cdot \frac{\rho _k -s}{\rho _k -\beta }} \right)} \\
&=\mathcal{H} \left( 0 \right)\prod\limits_{\rho _k } {\frac{\rho _k -\beta }{\rho
_k }} \cdot \prod\limits_{\rho _k } {\frac{\rho _k -s}{\rho _k -\beta }}\\
&=\mathcal{H} \left( 0 \right)\prod\limits_{\rho _k } {\left( {1-\frac{\beta }{\rho
_k }} \right)} \cdot \prod\limits_{\rho _k } {\frac{\rho _k -\beta -\left(
{s-\beta } \right)}{\rho _k -\beta }} \\
&=\mathcal{H} \left( 0 \right)\prod\limits_{\rho _k } {\left( {1-\frac{\beta }{\rho
_k }} \right)} \cdot \prod\limits_{\rho _k } {\left( {1-\frac{s-\beta }{\rho
_k -\beta }} \right)} . \\
 \end{aligned}
\end{equation}
Taking $s=\beta $ in (\ref{eq2}), we may get
\begin{equation}
\label{eq15}
\mathcal{H} \left( \beta \right)=\mathcal{H} \left( 0 \right)\prod\limits_{\rho _k }
{\left( {1-\frac{\beta }{\rho _k }} \right)} .
\end{equation}
By substituting (\ref{eq15}) into the last term of (\ref{eq14}), we carry out
\begin{equation}
\label{eq16}
\mathcal{H} \left( s \right)=\mathcal{H} \left( \beta \right)\prod\limits_{\rho _k }
{\left( {1-\frac{s-\beta }{\rho _k -\beta }} \right)} .
\end{equation}
Thus, we finish the proof of Lemma \ref{LE2}.
\end{proof}

\begin{corollary}
\label{COR1}
Assume $s\in \mathcal{C}$ , $\mathcal{H} \in \mathcal{J}$ and $\ell \in \mathbf{H}$. Then we have
\begin{equation}
\label{eq17}
\mathcal{H} \left( s \right)=\mathcal{H} \left( \ell \right)\prod\limits_{\rho _k }
{\left( {1-\frac{s-\ell }{\rho _k -\ell }} \right)}.
\end{equation}
\end{corollary}

\begin{proof}
Following the definition of $\mathcal{H} \in \mathcal{J}$, we have $\mathcal{H} \left( \ell
\right)\ne 0$ such that (\ref{eq17}) is valid if we make $\beta =\ell $ in
Lemma \ref{LE2}.
Thus, we finish the proof.
\end{proof}

\begin{corollary}
\label{COR2}
Suppose that $s\in \mathcal{C}$, $\mathcal{H} \in \mathcal{J}$ and $\ell \in \mathbf{H}$.
Then we have
\begin{equation}
\label{eq18}
\mathcal{H} \left( s \right)=\mathcal{H} \left( \ell \right)\prod\limits_{Im\left( {\rho
_k } \right)>0} {\left[ {1-\left( {\frac{s-\ell }{\rho _k -\ell }}
\right)^2} \right]}.
\end{equation}
\end{corollary}

\begin{proof}
By using Lemma \ref{LE1}, we have
\begin{equation}
\label{eq19}
\mathcal{H} \left( {\rho _k } \right)=\mathcal{H} \left( {2\ell -\rho _k } \right),
\end{equation}
which implies that
\begin{equation}
\label{eq20}
 \begin{aligned}
\mathcal{H} \left( s \right)&=\mathcal{H} \left( \ell \right)\prod\limits_{\rho _k }
{\left( {1-\frac{s-\ell }{\rho _k -\ell }} \right)} \\
&=\mathcal{H} \left( \ell \right)\prod\limits_{Im\left( {\rho _k } \right)>0}
{\left( {1-\frac{s-\ell }{\rho _k -\ell }} \right)\left[ {1-\frac{s-\ell
}{\left( {2\ell -\rho _k } \right)-\ell }} \right]} \\
&=\mathcal{H} \left( \ell \right)\prod\limits_{Im\left( {\rho _k } \right)>0}
{\left\{ {\left( {1-\frac{s-\ell }{\rho _k -\ell }} \right)\left(
{1-\frac{s-\ell }{\ell -\rho _k }} \right)} \right\}} \\
&=\mathcal{H} \left( \ell \right)\prod\limits_{Im\left( {\rho _k } \right)>0}
{\left\{ {\left( {1-\frac{s-\ell }{\rho _k -\ell }} \right)\left(
{1+\frac{s-\ell }{\rho _k -\ell }} \right)} \right\}} \\
&=\mathcal{H} \left( \ell \right)\prod\limits_{Im\left( {\rho _k } \right)>0}
{\left[ {1-\left( {\frac{s-\ell }{\rho _k -\ell }} \right)^2} \right]} . \\
 \end{aligned}
\end{equation}
\end{proof}

\begin{corollary}
\label{COR3}
If $s\in \mathcal{C}$, $\mathcal{H} \in \mathcal{J}$ and $\ell \in \mathbf{H}$, then there
exists any $\beta \in \mathcal{C}$ with $\beta \ne \rho _k $ such that
\begin{equation}
\label{eq21}
\mathcal{H} \left( s \right)=\mathcal{H} \left( \beta \right)\prod\limits_{Im\left( {\rho
_k } \right)>0} {\left[ {1-\frac{\left( {s-\ell } \right)^2-\left( {\beta
-\ell } \right)^2}{\left( {\rho _k -\ell } \right)^2-\left( {\beta -\ell }
\right)^2}} \right]}.
\end{equation}
\end{corollary}

\begin{proof}
By Lemmas \ref{LE1} and \ref{LE2}, we obtain
\begin{equation}
\label{eq22}
 \begin{aligned}
\mathcal{H} \left( s \right)&=\mathcal{H} \left( \beta \right)\prod\limits_{\rho _k }
{\left( {1-\frac{s-\beta }{\rho _k -\beta }} \right)}
 =\mathcal{H} \left( \beta \right)\prod\limits_{Im\left( {\rho _k } \right)>0}
{\left( {1-\frac{s-\beta }{\rho _k -\beta }} \right)\left( {1-\frac{s-\beta
}{2\ell -\rho _k -\beta }} \right)} \\
&=\mathcal{H} \left( \beta \right)\prod\limits_{Im\left( {\rho _k } \right)>0}
{\left[ {1-\frac{\left( {s-\ell } \right)-\left( {\beta -\ell }
\right)}{\left( {\rho _k -\ell } \right)-\left( {\beta -\ell } \right)}}
\right]\left[ {1-\frac{\left( {s-\ell } \right)+\left( {\ell -\beta }
\right)}{\left( {\ell -\rho _k } \right)+\left( {\ell -\beta } \right)}}
\right]} \\
&=\mathcal{H} \left( \beta \right)\prod\limits_{Im\left( {\rho _k } \right)>0}
{\left[ {\frac{\left( {\rho _k -\ell } \right)-\left( {s-\ell }
\right)}{\left( {\rho _k -\ell } \right)-\left( {\beta -\ell } \right)}}
\right]} \prod\limits_{Im\left( {\rho _k } \right)>0} {\left[ {\frac{\left(
{\ell -\rho _k } \right)-\left( {s-\ell } \right)}{\left( {\ell -\rho _k }
\right)+\left( {\ell -\beta } \right)}} \right]} , \\
\end{aligned}
\end{equation}
where
\begin{equation}
\label{eq23}
 \begin{aligned}
 \prod\limits_{Im\left( {\rho _k } \right)>0} {\left[ {1-\frac{\left(
{s-\ell } \right)-\left( {\beta -\ell } \right)}{\left( {\rho _k -\ell }
\right)-\left( {\beta -\ell } \right)}} \right]}
&=\prod\limits_{Im\left( {\rho _k } \right)>0} {\left[ {\frac{\left( {\rho
_k -\ell } \right)-\left( {\beta -\ell } \right)-\left( {s-\ell }
\right)+\left( {\beta -\ell } \right)}{\left( {\rho _k -\ell }
\right)-\left( {\beta -\ell } \right)}} \right]} \\
&=\prod\limits_{Im\left( {\rho _k } \right)>0} {\left[ {\frac{\left( {\rho
_k -\ell } \right)-\left( {s-\ell } \right)}{\left( {\rho _k -\ell }
\right)-\left( {\beta -\ell } \right)}} \right]} \\
 \end{aligned}
\end{equation}
and
\begin{equation}
\label{eq24}
 \begin{aligned}
 \prod\limits_{Im\left( {\rho _k } \right)>0} {\left[ {1-\frac{\left(
{s-\ell } \right)+\left( {\ell -\beta } \right)}{\left( {\ell -\rho _k }
\right)+\left( {\ell -\beta } \right)}} \right]}
&=\prod\limits_{Im\left( {\rho _k } \right)>0} {\left[ {\frac{\left( {\ell
-\rho _k } \right)+\left( {\ell -\beta } \right)-\left( {s-\ell }
\right)-\left( {\ell -\beta } \right)}{\left( {\ell -\rho _k }
\right)+\left( {\ell -\beta } \right)}} \right]} \\
&=\prod\limits_{Im\left( {\rho _k } \right)>0} {\left[ {\frac{\left( {\rho
_k -\ell } \right)+\left( {s-\ell } \right)}{\left( {\rho _k -\ell }
\right)+\left( {\beta -\ell } \right)}} \right]} . \\
 \end{aligned}
\end{equation}
To simplify (\ref{eq14}), we obtain
\begin{equation}
\label{eq25}
 \begin{aligned}
\mathcal{H} \left( s \right)&=\mathcal{H} \left( \beta \right)\prod\limits_{Im\left( {\rho
_k } \right)>0} {\left[ {\frac{\left( {\rho _k -\ell } \right)-\left(
{s-\ell } \right)}{\left( {\rho _k -\ell } \right)-\left( {\beta -\ell }
\right)}} \right]} \prod\limits_{Im\left( {\rho _k } \right)>0} {\left[
{\frac{\left( {\rho _k -\ell } \right)+\left( {s-\ell } \right)}{\left(
{\rho _k -\ell } \right)+\left( {\beta -\ell } \right)}} \right]} \\
&=\mathcal{H} \left( \beta \right)\prod\limits_{Im\left( {\rho _k } \right)>0}
{\left[ {\frac{\left( {\rho _k -\ell } \right)-\left( {s-\ell }
\right)}{\left( {\rho _k -\ell } \right)-\left( {\beta -\ell }
\right)}\frac{\left( {\rho _k -\ell } \right)+\left( {s-\ell }
\right)}{\left( {\rho _k -\ell } \right)+\left( {\beta -\ell } \right)}}
\right]} \\
&=\mathcal{H} \left( \beta \right)\prod\limits_{Im\left( {\rho _k } \right)>0}
{\left[ {\frac{\left( {\rho _k -\ell } \right)^2-\left( {s-\ell }
\right)^2}{\left( {\rho _k -\ell } \right)^2-\left( {\beta -\ell }
\right)^2}} \right]} \\
&=\mathcal{H} \left( \beta \right)\prod\limits_{Im\left( {\rho _k } \right)>0}
{\left[ {\frac{\left( {\rho _k -\ell } \right)^2-\left( {\beta -\ell }
\right)^2-\left( {s-\ell } \right)^2+\left( {\beta -\ell } \right)^2}{\left(
{\rho _k -\ell } \right)^2-\left( {\beta -\ell } \right)^2}} \right]} \\
&=\mathcal{H} \left( \beta \right)\prod\limits_{Im\left( {\rho _k } \right)>0}
{\left[ {1-\frac{\left( {s-\ell } \right)^2-\left( {\beta -\ell }
\right)^2}{\left( {\rho _k -\ell } \right)^2-\left( {\beta -\ell }
\right)^2}} \right]} \\
 \end{aligned}
\end{equation}
We then complete the proof.
\end{proof}

\subsection{The symmetric lines for $\mathcal{H} \left( s \right)$}

The symmetric function $\Lambda \left( s \right)$ is defined as
\begin{equation}
\label{eq26}
\Lambda \left( {s,\beta } \right)=:\frac{\mathcal{H} \left( s \right)}{\mathcal{H} \left(
\beta \right)}=\prod\limits_{Im\left( {\rho _k } \right)>0} {\left[
{1-\frac{\left( {s-\ell } \right)^2-\left( {\beta -\ell } \right)^2}{\left(
{\rho _k -\ell } \right)^2-\left( {\beta -\ell } \right)^2}} \right]} .
\end{equation}

\begin{corollary}
\label{COR4}
There exist the first symmetric line $s=\ell $ and the second symmetric line
$\beta =\ell $ for the symmetric function $\Lambda \left( s \right)$ for
$s\in \mathcal{C}$.
\end{corollary}

\begin{proof}
By using the definition of the symmetric function (\ref{eq26}), we reduce to the required result.
\end{proof}

\begin{corollary}
\label{COR5}
There exist the first symmetric line $s=\ell $ and the second symmetric line
$\beta =\ell $ for the entire function $\mathcal{H} \left( s \right)$ for $s\in \mathcal{C}
$.
\end{corollary}

\begin{proof}
Making use of (\ref{eq21}) in Corollary \ref{COR3}, we obtain the result.
\end{proof}

\begin{remark}
By (\ref{eq26}), we have
\begin{equation}
\label{eq27}
\Lambda \left( {s,\ell } \right)=\prod\limits_{Im\left( {\rho _k }
\right)>0} {\left[ {1-\left( {\frac{s-\ell }{\rho _k -\ell }} \right)^2}
\right]} ,
\end{equation}
and
\begin{equation}
\label{eq28}
\Lambda \left( {s,2\ell } \right)=\Lambda \left( {s,0}
\right)=\prod\limits_{Im\left( {\rho _k } \right)>0} {\left[ {1-\frac{\left(
{s-\ell } \right)^2-\ell ^2}{\left( {\rho _k -\ell } \right)^2-\ell ^2}}
\right]} ,
\end{equation}
such that
\begin{equation}
\label{eq29}
\mathcal{H} \left( s \right)=\mathcal{H} \left( 0 \right)\Lambda \left( {s,0} \right),
\end{equation}
\begin{equation}
\label{eq30}
\mathcal{H} \left( s \right)=\mathcal{H} \left( \ell \right)\Lambda \left( {s,\ell }
\right),
\end{equation}
and
\begin{equation}
\label{eq31}
\mathcal{H} \left( s \right)=\mathcal{H} \left( {2\ell } \right)\Lambda \left( {s,2\ell }
\right).
\end{equation}
Removing the effect of the second symmetric line $\beta =\ell $ and taking $\beta =\ell $ into
(\ref{eq13}) and (\ref{eq21}), we obtain (\ref{eq17}) and (\ref{eq18}), which can be applied to find
the critical line of the function $\mathcal{H} \left( s \right)$.
\end{remark}

\section{The proof of Theorem \ref{TH1}} {\label{sec:3}}

\subsection{A class of $\mathcal{H} \left( s
\right)$}

Since $\mathcal{H} \in \mathcal{J}$, $\mathcal{H} \left( s \right)$ can be
expressed as
\begin{equation}
\label{eq32}
\sum\limits_{m=0}^\infty {\Omega _m \left( {s-\ell } \right)^{2m}} =\mathcal{H}
\left( 0 \right)\prod\limits_{\rho _k } {\left( {1-\frac{s}{\rho _k }}
\right)} .
\end{equation}
If we combine (\ref{eq32}) and Lemma \ref{LE2}, then (\ref{eq32}) is equivalent to
\begin{equation}
\label{eq33}
\sum\limits_{m=0}^\infty {\Omega _m \left( {s-\ell } \right)^{2m}} =\mathcal{H}
\left( \beta \right)\prod\limits_{\rho _k } {\left( {1-\frac{s-\beta }{\rho
_k -\beta }} \right)} .
\end{equation}
Taking $\beta =\ell $ in (\ref{eq33}) implies that $\mathcal{H} \left( s \right)$ can be
rewritten as
\begin{equation}
\label{eq34}
\sum\limits_{m=0}^\infty {\Omega _m \left( {s-\ell } \right)^{2m}} =\mathcal{H}
\left( \ell \right)\prod\limits_{\rho _k } {\left( {1-\frac{s-\ell }{\rho _k
-\ell }} \right)} .
\end{equation}
Here, (\ref{eq34}) is obtained by (\ref{eq33}) when we remove the second symmetric line
$\beta =\ell $ in (\ref{eq33}).

By Corollary \ref{COR2} and (\ref{eq34}), we have
\begin{equation}
\label{eq35}
\mathcal{H} \left( \ell \right)\prod\limits_{\rho _k } {\left( {1-\frac{s-\ell
}{\rho _k -\ell }} \right)} =\mathcal{H} \left( \ell \right)\prod\limits_{Im\left(
{\rho _k } \right)>0} {\left[ {1-\left( {\frac{s-\ell }{\rho _k -\ell }}
\right)^2} \right]}
\end{equation}
such that
\begin{equation}
\label{eq36}
\sum\limits_{m=0}^\infty {\Omega _m \left( {s-\ell } \right)^{2m}} =\mathcal{H}
\left( \ell \right)\prod\limits_{Im\left( {\rho _k } \right)>0} {\left[
{1-\left( {\frac{s-\ell }{\rho _k -\ell }} \right)^2} \right]} .
\end{equation}

\begin{remark}
Clearly, the identity (\ref{eq36}) is also derived from Corollary \ref{COR3} when
we remove the second symmetric line $\beta =\ell $ in (\ref{eq21}) and is in
agreement with (\ref{eq6}).
\end{remark}

\subsection{A class of $\psi \left( \vartheta \right)$}

Taking $\vartheta =s-\ell $ in (\ref{eq36}) for $\vartheta \in \mathcal{C}$, we
get
\begin{equation}
\label{eq37}
\psi \left( \vartheta \right)=\sum\limits_{m=0}^\infty {\Omega _m \vartheta
^{2m}}
\end{equation}
and
\begin{equation}
\label{eq38}
\psi \left( \vartheta \right)=\mathcal{H} \left( \ell \right)\prod\limits_{Im\left(
{\rho _k } \right)>0} {\left[ {1-\frac{\vartheta ^2}{\left( {\rho _k -\ell }
\right)^2}} \right]} .
\end{equation}
Combining (\ref{eq37}) and (\ref{eq38}), we have the class of $\psi \left( \vartheta
\right)$, that is,
\begin{equation}
\label{eq39}
\sum\limits_{m=0}^\infty {\Omega _m \vartheta ^{2m}} =\mathcal{H} \left( \ell
\right)\prod\limits_{Im\left( {\rho _k } \right)>0} {\left[
{1-\frac{\vartheta ^2}{\left( {\rho _k -\ell } \right)^2}} \right]} .
\end{equation}
Obviously, $\psi \left( \vartheta \right)$ and $\mathcal{H} \left( s \right)$ are
entire functions of order $\nu =1$.

\subsection{A class of $\overline {\psi \left( \vartheta \right)}$}

It follows from (\ref{eq37}) that the complex conjugate $\overline {\psi \left(
\vartheta \right)} $ of the function $\psi \left( \vartheta \right)$ reads
\begin{equation}
\label{eq40}
\overline {\psi \left( \vartheta \right)} =\overline
{\sum\limits_{m=0}^\infty {\Omega _m \vartheta ^{2m}} }
=\sum\limits_{m=0}^\infty {\Omega _m \overline \vartheta ^{2m}}
\end{equation}
due to the fact
\begin{equation}
\label{eq41}
\Omega _m >0.
\end{equation}
From (\ref{eq40}) we show that
\begin{equation}
\label{eq42}
\overline {\psi \left( \vartheta \right)} =\psi \left( {\overline \vartheta
} \right).
\end{equation}
By using (\ref{eq39}) and (\ref{eq42}), we have
\begin{equation}
\label{eq43}
\begin{array}{l}
 \overline {\psi \left( \vartheta \right)} =\psi \left( {\overline \vartheta
} \right)=\sum\limits_{m=0}^\infty {\Omega _m \overline \vartheta ^{2m}}
 =\mathcal{H} \left( \ell \right)\prod\limits_{Im\left( {\rho _k } \right)>0}
{\left[ {1-\frac{\overline \vartheta ^2}{\left( {\rho _k -\ell } \right)^2}}
\right]}
 \end{array}
\end{equation}
because of $\mathcal{H} \left( \ell \right)=\Omega _0>0$.
Finding $\overline {\psi \left( \vartheta \right)} $ in (\ref{eq38}) implies that
\begin{equation}
\label{eq44}
 \begin{aligned}
 \overline {\psi \left( \vartheta \right)} &=\overline {\left\{ {\mathcal{H} \left(
\ell \right)\prod\limits_{Im\left( {\rho _k } \right)>0} {\left[
{1-\frac{\vartheta ^2}{\left( {\rho _k -\ell } \right)^2}} \right]} }
\right\}} =\mathcal{H} \left( \ell \right)\overline {\left\{
{\prod\limits_{Im\left( {\rho _k } \right)>0} {\left[ {1-\frac{\vartheta
^2}{\left( {\rho _k -\ell } \right)^2}} \right]} } \right\}} \\
&=\mathcal{H} \left( \ell \right)\prod\limits_{Im\left( {\rho _k } \right)>0}
{\left[ {1-\frac{\overline \vartheta ^2}{\left( {\overline {\rho _k } -\ell
} \right)^2}} \right]} . \\
 \end{aligned}
\end{equation}
Combining (\ref{eq43}) and (\ref{eq44}), we get
\begin{equation}
\label{eq45}
\begin{array}{l}
 \psi \left( {\overline \vartheta } \right)=\sum\limits_{m=0}^\infty {\Omega _m
\overline \vartheta ^{2m}} =\mathcal{H} \left( \ell
\right)\prod\limits_{Im\left( {\rho _k } \right)>0} {\left[
{1-\frac{\overline \vartheta ^2}{\left( {\rho _k -\ell } \right)^2}}
\right]}
 =\mathcal{H} \left( \ell \right)\prod\limits_{Im\left( {\rho _k } \right)>0}
{\left[ {1-\frac{\overline \vartheta ^2}{\left( {\overline {\rho _k } -\ell
} \right)^2}} \right]},
 \end{array}
\end{equation}
where $\vartheta \in \mathcal{C}$.

\subsection{Two products and convergence of $\psi \left( \vartheta \right)$}

By replacing $\overline \vartheta \in \mathcal{C}$ by $\vartheta \in \mathcal{C}$ in
(\ref{eq45}), we obtain
\begin{equation}
\label{eq46}
 \begin{aligned}
 \psi \left( \vartheta \right)&=\sum\limits_{m=0}^\infty {\Omega _m \vartheta
^{2m}} =\mathcal{H} \left( \ell \right)\prod\limits_{Im\left( {\rho _k } \right)>0}
{\left[ {1-\frac{\vartheta ^2}{\left( {\rho _k -\ell } \right)^2}} \right]}
=\mathcal{H} \left( \ell \right)\prod\limits_{Im\left( {\rho _k } \right)>0}
{\left[ {1-\frac{\vartheta ^2}{\left( {\overline {\rho _k } -\ell }
\right)^2}} \right]} .
 \end{aligned}
\end{equation}
Considering the fact that $\psi \left( \vartheta \right)$ is the even entire
function of order $\nu =1$, there exists any $r>0$ such that

\begin{equation}
\label{eq3.016}
\sum\limits_{\rho _k } {\left| {\rho _k -\ell } \right|^{-\left( {1+r}
\right)}}
\end{equation}
and
\begin{equation}
\label{eq3.017}
\sum\limits_{\rho _k } {\left| {\overline {\rho _k } -\ell }
\right|^{-\left( {1+r} \right)}}
\end{equation}
are convergent.

Taking $r=1$ in (\ref{eq3.016}) and (\ref{eq3.017}) implies that

\begin{equation}
\label{eq3.018}
\sum\limits_{\rho _k } {\left| {\rho _k -\ell } \right|^{-2}}
\end{equation}
and

\begin{equation}
\label{eq3.019}
\sum\limits_{\rho _k } {\left| {\overline {\rho _k } -\ell } \right|^{-2}}
\end{equation}
are absolutely convergent.

By (\ref{eq3.018}) and (\ref{eq3.019}), both
\begin{equation}
\label{eq47}
\sum\limits_{k=1}^\infty {\frac{1}{\left( {\rho _k -\ell } \right)^2}}
\end{equation}
and
\begin{equation}
\label{eq48}
\sum\limits_{k=1}^\infty {\frac{1}{\left( {\overline {\rho _k } -\ell }
\right)^2}} ,
\end{equation}
are convergent.

It follows from (\ref{eq46}), (\ref{eq47}) and (\ref{eq48}) that
\begin{equation}
\label{eq49}
\sum\limits_{k=1}^\infty {\frac{1}{\left( {\rho _k -\ell } \right)^2}}
=\sum\limits_{k=1}^\infty {\frac{1}{\left( {\overline {\rho _k } -\ell }
\right)^2}} .
\end{equation}

\subsection{A detailed proof of Theorem \ref{TH1}}

From (\ref{eq47}) we have
\begin{equation}
\label{eq50}
\left( {\rho _k -\ell } \right)^2=\left( {\overline {\rho _k } -\ell }
\right)^2,
\end{equation}
which leads to
\begin{equation}
\label{eq51}
\left[ {\left( {\rho _k -\ell } \right)+\left( {\overline {\rho _k } -\ell }
\right)} \right]\left[ {\left( {\rho _k -\ell } \right)-\left( {\overline
{\rho _k } -\ell } \right)} \right]=0.
\end{equation}
Since
\begin{equation}
\label{eq52}
\left( {\rho _k -\ell } \right)-\left( {\overline {\rho _k } -\ell }
\right)=\rho _k -\overline {\rho _k } =2iIm\left( {\rho _k } \right),
\end{equation}
we have from (\ref{eq51}) that
\begin{equation}
\label{eq53}
\left( {\rho _k -\ell } \right)+\left( {\overline {\rho _k } -\ell }
\right)=0.
\end{equation}
By (\ref{eq53}), we get
\begin{equation}
\label{eq54}
2Re\left( {\rho _k } \right)-2\ell =0
\end{equation}
or, alternatively,
\begin{equation}
\label{eq55}
Re\left( {\rho _k } \right)=\ell .
\end{equation}
We hence finish the proof of Theorem \ref{TH1}.

\section{An equivalent representation of Theorem \ref{TH1}} {\label{sec:4}}

Taking
\begin{equation}
\label{eq56}
Im\left( {\rho _k } \right)=\chi _k >0
\end{equation}
and using (\ref{eq45}) and (\ref{eq55}), we arrive at
\begin{equation}
\label{eq57}
 \begin{aligned}
 \psi \left( \vartheta \right)&=\sum\limits_{m=0}^\infty {\Omega _m \vartheta
^{2m}}
=\mathcal{H} \left( \ell \right)\prod\limits_{Im\left( {\rho _k } \right)>0}
{\left[ {1-\frac{\vartheta ^2}{\left( {\rho _k -\ell } \right)^2}} \right]}
=\mathcal{H} \left( \ell \right)\prod\limits_{\chi _k } {\left( {1+\frac{\vartheta
^2}{\chi _k^2 }} \right)} \\
&=\mathcal{H} \left( \ell \right)\prod\limits_{k=1}^\infty {\left(
{1+\frac{\vartheta ^2}{\chi _k^2 }} \right)} . \\
 \end{aligned}
\end{equation}
Clearly, all nontrivial zeros of $\psi \left( \vartheta \right)$ read $\pm
i\chi _k $, where $i=\sqrt {-1} $.

By the combination of (\ref{eq21}), (\ref{eq32}), (\ref{eq33}), (\ref{eq35}) and (\ref{eq55}), we have
\begin{equation}
\label{eq58}
 \begin{aligned}
\mathcal{H} \left( s \right)&=\sum\limits_{m=0}^\infty {\Omega _m \left( {s-\ell }
\right)^{2m}}
=\mathcal{H} \left( 0 \right)\prod\limits_{\rho _k } {\left(
{1-\frac{s}{\rho _k }} \right)}
=\mathcal{H} \left( \beta \right)\prod\limits_{\rho
_k } {\left( {1-\frac{s-\beta }{\rho _k -\beta }} \right)} \\
&=\mathcal{H} \left( \beta \right)\prod\limits_{Im\left( {\rho _k } \right)>0}
{\left[ {1-\frac{\left( {s-\ell } \right)^2-\left( {\beta -\ell }
\right)^2}{\left( {\rho _k -\ell } \right)^2-\left( {\beta -\ell }
\right)^2}} \right]} \\
&=\mathcal{H} \left( \ell \right)\prod\limits_{\rho _k } {\left( {1-\frac{s-\ell
}{\rho _k -\ell }} \right)} \\
&=\mathcal{H} \left( \ell \right)\prod\limits_{Im\left(
{\rho _k } \right)>0} {\left[ {1-\left( {\frac{s-\ell }{\rho _k -\ell }}
\right)^2} \right]} \\
&=\mathcal{H} \left( \ell \right)\prod\limits_{k=1}^\infty {\left( {1+\frac{\left(
{s-\ell } \right)^2}{\chi _k^2 }} \right)} . \\
 \end{aligned}
\end{equation}
Obviously, all nontrivial zeros of $\mathcal{H} \left( s \right)$ are written as
\begin{equation}
\label{eq59}
\rho _k =\ell \pm i\chi _k .
\end{equation}
By using (\ref{eq59}) and
\begin{equation}
\label{eq60}
\mathcal{H} \left( {\rho _k } \right)=\mathcal{H} \left( {2\ell -\rho _k } \right),
\end{equation}
we obtain
\begin{equation}
\label{eq61}
 \begin{aligned}
\mathcal{H} \left( s \right)&=\mathcal{H} \left( 0 \right)\prod\limits_{\rho _k } {\left(
{1-\frac{s}{\rho _k }} \right)}
=\mathcal{H} \left( 0 \right)\prod\limits_{Im\left(
{\rho _k } \right)>0} {\left( {1-\frac{s}{\rho _k }} \right)} \left(
{1-\frac{s}{2\ell -\rho _k }} \right) \\
&=\mathcal{H} \left( 0 \right)\prod\limits_{Im\left( {\rho _k } \right)>0} {\left\{
{\left( {1-\frac{s}{\rho _k }} \right)\left( {1-\frac{s}{\ell -i\chi _k }}
\right)} \right\}} \\
&=\mathcal{H} \left( 0 \right)\prod\limits_{Im\left( {\rho _k } \right)>0} {\left\{
{\left( {1-\frac{s}{\rho _k }} \right)\left( {1-\frac{s}{\ell -i\chi _k }}
\right)} \right\}} \\
&=\mathcal{H} \left( 0 \right)\prod\limits_{Im\left( {\rho _k } \right)>0} {\left\{
{\left( {1-\frac{s}{\rho _k }} \right)\left( {1-\frac{s}{\overline \rho _k
}} \right)} \right\}} . \\
 \end{aligned}
\end{equation}
It is obvious that $\rho _k $, $\overline \rho _k $, $2\ell -\rho _k $ and
$2\ell -\overline \rho _k $ are the nontrivial zeros of $\mathcal{H} \left( s
\right)$.

As a direct result of (\ref{eq58}), we have the following:

\begin{corollary}
\label{COR6}
Suppose $s\in \mathcal{C}$, $\mathcal{H} \in \mathcal{J}$ and $\ell \in \mathbf{H}$. Then
there exist the following equivalent representations:

($\mathcal{L}1$) All zeros of $\mathcal{H} \left( s \right)$ lie on the critical
line $Re\left( {\rho _k } \right)=\ell $.

($\mathcal{L}2$) There exists
\begin{equation}
\label{eq62}
\sum\limits_{m=0}^\infty {\Omega _m \left( {s-\ell } \right)^{2m}} =\mathcal{H}
\left( 0 \right)\prod\limits_{\rho _k } {\left( {1-\frac{s}{\rho _k }}
\right)} .
\end{equation}

($\mathcal{L}3$) There exists
\begin{equation}
\label{eq63}
\sum\limits_{m=0}^\infty {\Omega _m \left( {s-\ell } \right)^{2m}} =\mathcal{H}
\left( \ell \right)\prod\limits_{\rho _k } {\left( {1-\frac{s-\ell }{\rho _k
-\ell }} \right)} .
\end{equation}

($\mathcal{L}4$) There exists
\begin{equation}
\label{eq64}
\sum\limits_{m=0}^\infty {\Omega _m \left( {s-\ell } \right)^{2m}} =\mathcal{H}
\left( \ell \right)\prod\limits_{Im\left( {\rho _k } \right)>0} {\left[
{1-\left( {\frac{s-\ell }{\rho _k -\ell }} \right)^2} \right]} .
\end{equation}

($\mathcal{L}5$) There exists
\begin{equation}
\label{eq65}
\sum\limits_{m=0}^\infty {\Omega _m \left( {s-\ell } \right)^{2m}} =\mathcal{H}
\left( \ell \right)\prod\limits_{k=1}^\infty {\left( {1+\frac{\left( {s-\ell
} \right)^2}{\chi _k^2 }} \right)} .
\end{equation}

($\mathcal{L}6$) There exists
\begin{equation}
\label{eq66}
\sum\limits_{m=0}^\infty {\Omega _m \left( {s-\ell } \right)^{2m}} =\mathcal{H}
\left( \beta \right)\prod\limits_{\rho _k } {\left( {1-\frac{s-\beta }{\rho
_k -\beta }} \right)} .
\end{equation}

($\mathcal{L}7$) There exists
\begin{equation}
\label{eq67}
\sum\limits_{m=0}^\infty {\Omega _m \left( {s-\ell } \right)^{2m}} =\mathcal{H}
\left( \beta \right)\prod\limits_{Im\left( {\rho _k } \right)>0} {\left[
{1-\frac{\left( {s-\ell } \right)^2-\left( {\beta -\ell } \right)^2}{\left(
{\rho _k -\ell } \right)^2-\left( {\beta -\ell } \right)^2}} \right]} .
\end{equation}
\end{corollary}

\begin{remark}
As a matter of fact, the equivalent relationship between ($\mathcal{L}1$) and ($\mathcal{L}2$)
was adopted in \cite{20} if $\mathcal{H} \left( s \right)$  is considered as the Riemann xi function.
The result considered in \cite {21} is ($\mathcal{L}4$) at $\ell=0$. In Corollary \ref{COR6},
we adopt ($\mathcal{L}7$) to reduce to ($\mathcal{L}4$) because this process removes
the influence of the the second symmetric line
$\beta =\ell $ for the entire function $\mathcal{H} \left( s \right)$. There always exists
an entire function $\psi \left( \vartheta \right)$ of order $\nu =1$, i.e.,
\begin{equation}
\label{eqA1}
 \begin{aligned}
\psi \left( \vartheta \right)=\sum\limits_{m=0}^\infty {\Omega _m \vartheta
^{2m}}
 \end{aligned}
\end{equation}
such that
\begin{equation}
\label{eqA2}
 \begin{aligned}
 \psi \left( \vartheta \right)=\mathcal{H} \left( \ell \right)\prod\limits_{Im\left( {\rho _k } \right)>0}
{\left[ {1-\frac{\vartheta ^2}{\left( {\rho _k -\ell } \right)^2}} \right]}
=\mathcal{H} \left( \ell \right)\prod\limits_{Im\left( {\rho _k } \right)>0}
{\left[ {1-\frac{\vartheta ^2}{\left( {\overline {\rho _k } -\ell }
\right)^2}} \right]}.
 \end{aligned}
\end{equation}
Since $\psi \left( \vartheta \right)$ is an even entire function of order $\nu =1$, this implies that
\begin{equation}
\label{eqA3}
\sum\limits_{k=1}^\infty {\frac{1}{\left( {\rho _k -\ell } \right)^2}}
=\sum\limits_{k=1}^\infty {\frac{1}{\left( {\overline {\rho _k } -\ell }
\right)^2}}
\end{equation}
or, alternatively,
\begin{equation}
\label{eqA4}
Re\left( {\rho _k } \right)=\ell .
\end{equation}
If we allow to take the value $\ell=0$ in (\ref{eqA2}), (\ref{eqA3}) and (\ref{eqA4}),
this is the key work presented in \cite {21}.
Thus, Corollary \ref{COR6} is an equivalent representation theorem
for the critical line for the entire function $\mathcal{H} \left( s \right)$
considered in the class $\mathcal{J}$.
\end{remark}
\section{A typical application associated with the work of Euler} {\label{sec:5}}

Suppose the special hyperbolic cosine function  is represented in the form
\begin{equation}
\label{eq68}
F\left( s \right)=cosh\left( {s-6} \right),
\end{equation}
where $s\in \mathcal{C}$.

From (\ref{eq68}) we obtain the series representation as follows:
\begin{equation}
\label{eq69}
F\left( s \right)=\sum\limits_{m=0}^\infty {\frac{\left( {s-6}
\right)^{2m}}{\left( {2m} \right)!}} .
\end{equation}
It follows from (\ref{eq69}) that
\begin{equation}
\label{eq70}
F\left( s \right)=F\left( {12-s} \right).
\end{equation}
Let us recall that
\begin{equation}
\label{eq71}
\cos \left( \vartheta \right)=\sum\limits_{m=0}^\infty {\frac{\left( {-1}
\right)^m}{\left( {2m} \right)!}\vartheta ^{2m}}
\end{equation}
and
\begin{equation}
\label{eq72}
\cosh \left( \vartheta \right)=\sum\limits_{m=0}^\infty {\frac{\vartheta
^{2m}}{\left( {2m} \right)!}} .
\end{equation}
Since $\cos \left( \vartheta \right)$ is an even function of order $\nu
=1$, there exist
\begin{equation}
\label{eq73}
\cosh \left( \vartheta \right)=\prod\limits_{\lambda _k } {\left(
{1-\frac{s}{\lambda _k }} \right)}
\end{equation}
and
\begin{equation}
\label{eq74}
\cos \left( \vartheta \right)=\prod\limits_{\lambda _k } {\left(
{1+\frac{s}{i\lambda _k }} \right)} ,
\end{equation}
where $\lambda _k \in \mathcal{C}$ run over the zeros of $\cosh \left( \vartheta
\right)$.

By combining (\ref{eq69}), (\ref{eq72}) and (\ref{eq74}), we suggest
\begin{equation}
\label{eq75}
F\left( s \right)=\cosh \left( {s-6} \right)=\sum\limits_{m=0}^\infty
{\frac{\left( {s-6} \right)^{2m}}{\left( {2m} \right)!}}
=\prod\limits_{\lambda _k } {\left( {1-\frac{s-6}{\lambda _k }} \right)} .
\end{equation}
From (\ref{eq75}) we reduce to
\begin{equation}
\label{eq76}
\phi _k =\lambda _k +6.
\end{equation}
By using (\ref{eq76}), we write (\ref{eq75}) as
\begin{equation}
\label{eq77}
F\left( s \right)=\cosh \left( {s-6} \right)=\sum\limits_{m=0}^\infty
{\frac{\left( {s-6} \right)^{2m}}{\left( {2m} \right)!}} =\prod\limits_{\phi
_k } {\left( {1-\frac{s-6}{\phi _k -6}} \right)} ,
\end{equation}
which leads to
\begin{equation}
\label{eq78}
F\left( 6 \right)=1.
\end{equation}
Combining (\ref{eq77}) and (\ref{eq78}), $F\left( s \right)=\cosh \left( {s-6} \right)$
can be represented as
\begin{equation}
\label{eq79}
\sum\limits_{m=0}^\infty {\frac{\left( {s-6} \right)^{2m}}{\left( {2m}
\right)!}} =F\left( 6 \right)\prod\limits_{\phi _k } {\left(
{1-\frac{s-6}{\phi _k -6}} \right)} .
\end{equation}
Let us recall that $\cos \left( \vartheta \right)$ is an even function of
order $\nu=1$. Then, by (\ref{eq73}), there exists any $\hbar >0$ such that
$\sum\limits_{\lambda _k } {\left| {\lambda _k } \right|^{-\left( {1+\hbar }
\right)}}$ is convergent. This implies that we have from (\ref{eq76}) that
\begin{equation}
\label{eq80}
\sum\limits_{\lambda _k } {\left| {\lambda _k } \right|^{-\left( {1+\hbar }
\right)}} =\sum\limits_{\phi _k } {\left| {\phi _k -6} \right|^{-\left(
{1+\hbar } \right)}}
\end{equation}
is convergent.

This implies that $F\in \mathcal{J}$ and that by (A3) in Corollary \ref{COR6} we
have
\begin{equation}
\label{eq81}
Re\left( {\phi _k } \right)=6.
\end{equation}
From (\ref{eq76}) and (\ref{eq81}) we deduce that
\begin{equation}
\label{eq82}
Re\left( {\lambda _k } \right)=0.
\end{equation}
Let $Im\left( {\lambda _k } \right)=\mu _k >0$. Then we have from (\ref{eq74}) and
(\ref{eq82}) that
\begin{equation}
\label{eq83}
 \begin{aligned}
 \cos \left( \vartheta \right)&=\prod\limits_{\lambda _k } {\left(
{1+\frac{\vartheta }{i\lambda _k }} \right)}
=\prod\limits_{Im\left(
{\lambda _k } \right)>0} {\left( {1+\frac{\vartheta }{i\lambda _k }}
\right)} \left( {1-\frac{\vartheta }{i\lambda _k }} \right) \\
&=\prod\limits_{k=1}^\infty {\left[ {\left( {1+\frac{\vartheta }{i\left(
{i\mu _k } \right)}} \right)\left( {1-\frac{\vartheta }{i\left( {i\mu _k }
\right)}} \right)} \right]} \\
&=\prod\limits_{k=1}^\infty {\left[ {\left(
{1-\frac{\vartheta }{\mu _k }} \right)\left( {1+\frac{\vartheta }{\mu _k }}
\right)} \right]} \\
&=\prod\limits_{k=1}^\infty {\left( {1-\frac{\vartheta ^2}{\mu _k^2 }}
\right)} . \\
 \end{aligned}
\end{equation}
Adopting (\ref{eq83}), we obtain $\mu _k^2 >0$, which agrees with the result of
Euler, i.e.,
\begin{equation}
\label{eq84}
\left( {2k-1} \right)^2\pi ^2/4>0
\end{equation}
since there exists (see \cite{22}; also see \cite{23}, p.114)
\begin{equation}
\label{eq85}
\cos \left( \vartheta \right)=\sum\limits_{m=0}^\infty {\frac{\left( {-1}
\right)^m\vartheta ^{2m}}{\left( {2m} \right)!}} =\prod\limits_{k=1}^\infty
{\left[ {1-\frac{4\vartheta ^2}{\left( {2k-1} \right)^2\pi ^2}} \right]} .
\end{equation}
From (\ref{eq76}) and (\ref{eq79}) we conclude that
\begin{equation}
\label{eq86}
\phi _k =6\pm i\mu _k .
\end{equation}
Considering $Im\left( {\lambda _k } \right)=\mu _k >0$ and using (\ref{eq86}),
the identity (\ref{eq79}) yields that
\begin{equation}
\label{eq87}
 \begin{aligned}
 F\left( s \right)&=\sum\limits_{m=0}^\infty {\frac{\left( {s-6}
\right)^{2m}}{\left( {2m} \right)!}} \\
&=F\left( 6 \right)\prod\limits_{\phi _k
} {\left( {1-\frac{s-6}{\phi _k -6}} \right)} \\
&=F\left( 6 \right)\prod\limits_{Im\left( {\phi _k } \right)>0} {\left(
{1-\frac{s-6}{\phi _k -6}} \right)\left( {1-\frac{s-6}{6-\phi _k }} \right)}
\\
&=F\left( 6 \right)\prod\limits_{Im\left( {\phi _k } \right)>0} {\left(
{1-\frac{s-6}{\phi _k -6}} \right)\left( {1+\frac{s-6}{\phi _k -6}} \right)}\\
&=F\left( 6 \right)\prod\limits_{Im\left( {\phi _k } \right)>0} {\left[
{1-\left( {\frac{s-6}{\phi _k -6}} \right)^2} \right]} \\
&=F\left( 6 \right)\prod\limits_{k=1}^\infty {\left[ {1-\left(
{\frac{s-6}{i\mu _k }} \right)^2} \right]} . \\
 \end{aligned}
\end{equation}
To simplify (\ref{eq87}), we obtain
\begin{equation}
\label{eq88}
F\left( s \right)=\sum\limits_{m=0}^\infty {\frac{\left( {s-6}
\right)^{2m}}{\left( {2m} \right)!}} =F\left( 6
\right)\prod\limits_{k=1}^\infty {\left[ {1+\frac{\left( {s-6}
\right)^2}{\mu _k^2 }} \right]} .
\end{equation}
Comparing between (\ref{eq85}) and (\ref{eq88}), we present
\begin{equation}
\label{eq89}
\mu _k =\frac{\left( {2k-1} \right)\pi }{2}.
\end{equation}
Thus,
\begin{equation}
\label{eq90}
F\left( s \right)=\sum\limits_{m=0}^\infty {\frac{\left( {s-6}
\right)^{2m}}{\left( {2m} \right)!}} =F\left( 6
\right)\prod\limits_{k=1}^\infty {\left[ {1+\frac{4\left( {s-6}
\right)^2}{\left( {2k-1} \right)^2\pi ^2}} \right]}
\end{equation}
and
\begin{equation}
\label{eq91}
F\left( s \right)=\sum\limits_{m=0}^\infty {\frac{\left( {s-6}
\right)^{2m}}{\left( {2m} \right)!}} =F\left( 6
\right)\prod\limits_{k=1}^\infty {\left[ {1+\frac{4\left( {s-6}
\right)^2}{\left( {2k-1} \right)^2\pi ^2}} \right]}
\end{equation}
With (\ref{eq88}) and (\ref{eq89}), we may get
\begin{equation}
\label{eq92}
\sum\limits_{k=1}^\infty {\mu _k^{-2} } =\sum\limits_{k=1}^\infty {\left|
{\frac{\left( {2k-1} \right)\pi }{2}} \right|^{-2}} =2\pi
^{-2}\sum\limits_{k=1}^\infty {\left( {k-\frac{1}{2}} \right)^{-2}} .
\end{equation}
Define the Hurwitz zeta function $\zeta \left( {\eta ,s} \right)$ by (see \cite{24}, p.607)
\begin{equation}
\label{eq93}
\zeta \left( {\eta ,s} \right)=\sum\limits_{m=0}^\infty {\left( {m+\eta }
\right)^{-s}} .
\end{equation}
In view of (\ref{eq92}) and (\ref{eq93}), we have
\begin{equation}
\label{eq94}
\sum\limits_{k=1}^\infty {\left( {k-\frac{1}{2}} \right)^{-2}} =\zeta \left(
{-\frac{3}{2},2} \right)-\frac{4}{9}
\end{equation}
such that
\begin{equation}
\label{eq95}
\sum\limits_{k=1}^\infty {\mu _k^{-2} } =2\pi ^{-2}\left[ {\zeta \left(
{-\frac{3}{2},2} \right)-\frac{4}{9}} \right].
\end{equation}
By using (\ref{eq86}) and ($\mathcal{L}7$) in Corollary \ref{COR6}, we may find that
\begin{equation}
\label{eq96}
 \begin{aligned}
 F\left( s \right)&=\sum\limits_{m=0}^\infty {\frac{\left( {s-6}
\right)^{2m}}{\left( {2m} \right)!}} \\
&=\Phi \left( \beta\right)\prod\limits_{Im\left( {\phi _k } \right)>0} {\left[ {1-\frac{\left(
{s-6} \right)^2-\left( {\beta -6} \right)^2}{\left( {\phi _k -6}
\right)^2-\left( {\beta -6} \right)^2}} \right]} \\
&=\Phi \left( \beta \right)\prod\limits_{k=1}^\infty {\left[ {1+\frac{\left(
{s-6} \right)^2-\left( {\beta -6} \right)^2}{\mu _k^2 +\left( {\beta -6}
\right)^2}} \right]}\\
&=\Phi \left( \beta \right)\prod\limits_{k=1}^\infty {\left[ {1+\frac{\left(
{s-6} \right)^2-\left( {\beta -6} \right)^2}{{\frac{\left( {2k-1} \right)^2\pi^2}{4}}+\left( {\beta -6}
\right)^2}} \right]} .\\
 \end{aligned}
\end{equation}
It follows from (\ref{eq96}) that $F\left( s \right)$ has the first symmetric line
$s=6$ and the second symmetric line $\beta =6$.

As a direct result of Corollary \ref{COR3}, we obtain the followings:

\begin{corollary}
\label{COR7}
If $s\in \mathcal{C}$ and $\ell \in \mathbf{H}$, then there
exists any $\beta \in \mathcal{C}$ with $\beta  \ne \ell  + i\left( {2k - 1} \right) \pi /2$ such that
\begin{equation}
\label{eq97}
\cosh \left( {s - \ell } \right) = \cosh \left( \beta  \right)\prod\limits_{k = 1}^\infty  {\left[ {1 + \frac{{\left( {s - \ell } \right)^2  - \left( {\beta  - \ell } \right)^2 }}{{\frac{{\left( {2k - 1} \right)^2 \pi ^2 }}{4} + \left( {\beta  - \ell } \right)^2 }}} \right]}.
\end{equation}
\end{corollary}

\begin{corollary}
\label{COR8}
If $s\in \mathcal{C}$ and $\ell \in \mathbf{H}$, then there
exists any $\beta \in \mathcal{C}$ with $\beta  \ne \ell  + ik\pi$ such that
\begin{equation}
\label{eq98}
{\mathop{\rm sinch}\nolimits} \left( {s - \ell } \right) = {\mathop{\rm sinch}\nolimits} \left( \beta  \right)\prod\limits_{k = 1}^\infty  {\left[ {1 + \frac{{\left( {s - \ell } \right)^2  - \left( {\beta  - \ell } \right)^2 }}{{k^2 \pi ^2  + \left( {\beta  - \ell } \right)^2 }}} \right]},
\end{equation}
where
\[
{\mathop{\rm sinch}\nolimits}\left( s \right) = \sum\limits_{m = 0}^\infty  {\frac{{s^{2n} }}{{\left( {2n + 1} \right)!}}}.
\]
\end{corollary}

\begin{corollary}
\label{COR9}
If $s\in \mathcal{C}$ and $\ell \in \mathbf{H}$, then there
exists any $\beta \in \mathcal{C}$ with $\beta  \ne \ell  + \left( {2k - 1} \right) \pi /2$ such that
\begin{equation}
\label{eq99}
\cos \left( {s - \ell } \right) = \cos \left( \beta  \right)\prod\limits_{k = 1}^\infty  {\left[ {1 - \frac{{\left( {s - \ell } \right)^2  - \left( {\beta  - \ell } \right)^2 }}{{\frac{{\left( {2k - 1} \right)^2 \pi ^2 }}{4} - \left( {\beta  - \ell } \right)^2 }}} \right]}.
\end{equation}
\end{corollary}

\begin{corollary}
\label{COR10}
If $s\in \mathcal{C}$ and $\ell \in \mathbf{H}$, then there
exists any $\beta \in \mathcal{C}$ with $\beta  \ne \ell  + k\pi$ such that
\begin{equation}
\label{eq100}
{\mathop{\rm sinc}\nolimits} \left( {s - \ell } \right) = {\mathop{\rm sinc}\nolimits} \left( \beta  \right)\prod\limits_{k = 1}^\infty  {\left[ {1 - \frac{{\left( {s - \ell } \right)^2  - \left( {\beta  - \ell } \right)^2 }}{{k^2 \pi ^2  - \left( {\beta  - \ell } \right)^2 }}} \right]} ,
\end{equation}
where ${\mathop{\rm sinc}\nolimits} \left( s \right)$ is the sinc function, denoted by \cite{25}
\[
{\mathop{\rm sinc}\nolimits} \left( s \right) = \sum\limits_{m = 0}^\infty  {\frac{{\left( { - 1} \right)^n s^{2n} }}{{\left( {2n + 1} \right)!}}}.
\]
\end{corollary}

\begin{remark}
The identity (\ref{eq96}) is a special case of (\ref{eq97}) at the point $\ell=6$.
Taking $\ell=0$ into (\ref{eq97}), (\ref{eq98}), (\ref{eq99}) and (\ref{eq100}),
we may have the followings:
\begin{equation}
\label{eq101}
\cosh \left( {s } \right) = \cosh \left( \beta  \right)\prod\limits_{k = 1}^\infty  {\left[ {1 + \frac{{s^2  - \beta^2 }}{{\frac{{\left( {2k - 1} \right)^2 \pi ^2 }}{4} + \beta^2 }}} \right]},
\end{equation}

\begin{equation}
\label{eq102}
{\mathop{\rm sinch}\nolimits} \left( {s} \right) = {\mathop{\rm sinch}\nolimits} \left( \beta  \right)\prod\limits_{k = 1}^\infty  {\left[ {1 + \frac{{s^2  - \beta^2 }}{{k^2 \pi ^2  + \beta^2 }}} \right]},
\end{equation}

\begin{equation}
\label{eq103}
\cos \left( {s} \right) = \cos \left( \beta  \right)\prod\limits_{k = 1}^\infty  {\left[ {1 - \frac{{s^2  - \beta^2 }}{{\frac{{\left( {2k - 1} \right)^2 \pi ^2 }}{4} - \beta^2 }}} \right]}
\end{equation}
and
\begin{equation}
\label{eq104}
{\mathop{\rm sinc}\nolimits} \left( {s} \right) = {\mathop{\rm sinc}\nolimits} \left( \beta  \right)\prod\limits_{k = 1}^\infty  {\left[ {1 - \frac{{s^2  - \beta^2 }}{{k^2 \pi ^2  - \beta^2 }}} \right]}.
\end{equation}

Putting $\beta=0$ into (\ref{eq101}), (\ref{eq102}), (\ref{eq103}) and (\ref{eq104}),
we can obtain
\begin{equation}
\label{eq105}
\cosh \left( {s } \right) = \cosh \left( 0  \right)\prod\limits_{k = 1}^\infty  {\left[ {1 + \frac{{s^2}}{{\frac{{\left( {2k - 1} \right)^2 \pi ^2 }}{4}}}} \right]}=\prod\limits_{k = 1}^\infty  {\left[ {1 + \frac{{s^2}}{{\frac{{\left( {2k - 1} \right)^2 \pi ^2 }}{4}}}} \right]},
\end{equation}

\begin{equation}
\label{eq106}
{\mathop{\rm sinch}\nolimits} \left( {s} \right) = {\mathop{\rm sinch}\nolimits} \left( 0 \right)\prod\limits_{k = 1}^\infty  {\left( {1 + \frac{{s^2}}{{k^2 \pi ^2}}} \right)}= \prod\limits_{k = 1}^\infty  {\left( {1 + \frac{{s^2}}{{k^2 \pi ^2}}} \right)},
\end{equation}

\begin{equation}
\label{eq107}
\cos \left( {s} \right) = \cos \left(0  \right)\prod\limits_{k = 1}^\infty  {\left[ {1 - \frac{{s^2}}{{\frac{{\left( {2k - 1} \right)^2 \pi ^2 }}{4}}}} \right]}= \prod\limits_{k = 1}^\infty  {\left[ {1 - \frac{{s^2}}{{\frac{{\left( {2k - 1} \right)^2 \pi ^2 }}{4}}}} \right]}
\end{equation}
and
\begin{equation}
\label{eq108}
{\mathop{\rm sinc}\nolimits} \left( {s} \right) = {\mathop{\rm sinc}\nolimits} \left(0  \right)\prod\limits_{k = 1}^\infty  {\left({1 - \frac{{s^2}}{{k^2 \pi ^2}}} \right)}=\prod\limits_{k = 1}^\infty  {\left({1 - \frac{{s^2}}{{k^2 \pi ^2}}} \right)}.
\end{equation}

Obviously, (\ref{eq105}) and (\ref{eq107}) are the results of Euler (see \cite{22}; also see \cite{24}, p.126 and p.118), and
(\ref{eq108}) is the result reported in \cite{25}.

By using (\ref{eq102}), (\ref{eq104}), (\ref{eq105}) and (\ref{eq107}), we may carry out
\begin{equation}
\label{eq109}
\begin{array}{l}
{\mathop{\rm sinh}\nolimits} \left( {s} \right)
={s}\times{\mathop{\rm sinch}\nolimits} \left( {s} \right)
= {\mathop{\rm sinh}\nolimits} \left( \beta  \right)s\prod\limits_{k = 1}^\infty  {\left[ {1 + \frac{{s^2  - \beta^2 }}{{k^2 \pi ^2  + \beta^2 }}} \right]},
\end{array}
\end{equation}

\begin{equation}
\label{eq110}
\begin{array}{l}
{\mathop{\rm sin}\nolimits} \left( {s} \right)
=s\times{\mathop{\rm sinc}\nolimits} \left( {s} \right)
={\mathop{\rm sinc}\nolimits} \left( \beta  \right) s \prod\limits_{k = 1}^\infty  {\left[ {1 - \frac{{s^2  - \beta^2 }}{{k^2 \pi ^2  - \beta^2 }}} \right]},
\end{array}
\end{equation}

\begin{equation}
\label{eq111}
\begin{array}{l}
{\mathop{\rm sinh}\nolimits} \left( {s} \right)
= {s} \prod\limits_{k = 1}^\infty  {\left( {1 + \frac{{s^2}}{{k^2 \pi ^2}}} \right)}
\end{array}
\end{equation}
and
\begin{equation}
\label{eq112}
\begin{array}{l}
{\mathop{\rm sin}\nolimits} \left( {s} \right)
=s\prod\limits_{k = 1}^\infty  {\left({1 - \frac{{s^2}}{{k^2 \pi ^2}}} \right)}.
\end{array}
\end{equation}
Clearly, (\ref{eq111}) and (\ref{eq112}) are the results of Euler (see \cite{22}; also see \cite{24}, p.126 and p.118).
\end{remark}
As a similar way of (\ref{eq92}), we have
\begin{equation}
\label{eq113}
\sum\limits_{k=1}^\infty {\left( {k\pi} \right)^{-2}} ={\pi^{-2}}\zeta \left(
{2} \right),
\end{equation}
where $\zeta \left({s} \right)$ is the Riemann zeta function, denoted as (see \cite{23}, p.151)
\[
\zeta \left({s} \right)=\sum\limits_{k=1}^\infty {k^{-s}}.
\]
It is obvious that (\ref{eq113}) holds because ${\mathop{\rm sinch}\nolimits} \left( {s} \right)$ is en even function of order one.

As a direct result, from Corollaries \ref{COR7} and \ref{COR8} we obtain the following:
\begin{corollary}
\label{COR11}
Assume the notations of Corollaries \ref{COR7} and \ref{COR8}.
Then $\cosh \left( {s - \ell } \right)$  and ${\mathop{\rm sinch}\nolimits} \left( {s - \ell } \right)$  belong to the class $\mathcal{J}$.
\end{corollary}
\section{Conclusion and further remarks} {\label{sec:6}}

In the present work we proposed a new class of the entire function of order
one with the real positive coefficients and infinity of complex zeros. We
suggested a sufficient condition for the same critical line for its complex
zeros. We presented the equivalent representation theorem for it. We also gave
a typical example for the special hyperbolic cosine function, whose result
is in accord with the result of Euler. The obtained result is proposed as a
new mathematical tool to obtain the critical line of the class of the entire
function of order one.


\begin{thebibliography}{120}

\bibitem{1}
B. Y. Levin, Distribution of zeros of entire functions, Vol. 150, American Mathematical Society, 1980.

\bibitem{2}
I. Markushevich, Entire functions, Elsevier, 2014.

\bibitem{3}
E. Laguerre, Sur les fonctions du genre z\'{e}ro et du genre un, Comptes rendus de l'Acad'emie des Sciences Paris, 95 (1882), 828-831.

\bibitem{4}
G. P\'{o}lya, \"{U}ber Ann\"{a}herung durch Polynome mit lauter reellen Wurzeln, Rendiconti del Circolo Matematico di Palermo, 36 (1913) (2), 279-295.

\bibitem{5}
I. Wagner, On a new class of Laguerre-P\'{o}lya type functions with applications in number theory, Pacific Journal of Mathematics, 320 (2022) (2), 177-192.

\bibitem{6}
\'{A}. Baricz and S. Singh, Zeros of some special entire functions, Proceedings of the American Mathematical Society, 146 (2018) (6), 2207-2216.

\bibitem{7}
G. Csordas and A. Vishnyakova, The generalized Laguerre inequalities and functions in the Laguerre-P\'{o}lya class, Central European Journal of Mathematics, 11 (2013) (10), 1643-1650.

\bibitem{8}
G. Csordas and D. K. Dimitrov, Conjectures and theorems in the theory of entire functions, Numerical Algorithms, 25 (2000) (2), 109-122.

\bibitem{9}
Y. O. Kim, Critical points of real entire functions and a conjecture of P\'{o}lya, Proceedings of the American Mathematical Society, 124 (1996) (4), 819-830.

\bibitem{10}
G. Csordas, R. S. Varga and I. Vincze, Jensen polynomials with applications to the Riemann $\xi $-function, Journal of Mathematical Analysis and Applications, 153 (1990) (2), 112-135.

\bibitem{11}
G. P\'{o}lya, Sur une question concernant les fonctions entieres, Comptes Rendus de l'Acad\'{e}mie des Sciences, Paris, 158 (1914), 330-333.

\bibitem{12}
S. Hellerstein and J. Williamson, Derivatives of entire functions and a question of P\'{o}lya, Transactions of the American Mathematical Society, 227 (1977), 227-249.

\bibitem{13}
N. G. De Bruijn, The roots of trigonometric integrals, Duke Mathematical Journal, 17 (1950) (4), 197-226.

\bibitem{14}
C. M. Newman, Fourier transforms with only real zeros, Proceedings of the American Mathematical Society, 61(1976) (3), 245-251.

\bibitem{15}
S. Karlin, Total positivity, Vol. 1, Stanford University Press, 1968.

\bibitem{16}
L. C. Shen, On the zeros of successive derivatives of even Laguerre-P\'{o}lya functions, Transactions of the American Mathematical Society, 298 (1986) (3), 643-652.

\bibitem{17}
D. Su\'{a}rez, A generalization of the Laguerre-P\'{o}lya class of entire functions, Journal of approximation theory, 101 (1999) (2), 37-48.

\bibitem{18}
D. K. Dimitrov and Y. B. Cheikh, Laguerre polynomials as Jensen polynomials of Laguerre-P\'{o}lya entire functions, Journal of Computational and Applied Mathematics, 233 (2009) (4), 703-707.

\bibitem{19}
A. Bohdanov and A. Vishnyakova, On the conditions for entire functions related to the partial theta function to belong to the Laguerre-P\'{o}lya class, Journal of Mathematical Analysis and Applications, 434 (2016) (3), 1740-1752.

\bibitem{20}
X. J. Yang, All nontrivial zeros for the Riemann zeta function are on the critical line $\Re (s) = 1/2$, arXiv: 1811.02418v15.

\bibitem{21}
X. J. Yang, On all real zeros for a class of the even entire function, arXiv: 2107.04005v2.

\bibitem{22}
L. Euler, Introductio in analysin infinitorum, Apud Marcum-Michaelem Bousquet, 1748.

\bibitem{23}
E. C. Titchmarsh, The theory of functions, Oxford University Press, 1939.

\bibitem{24}
F. W. Olver, D. W. Lozier, R. F. Boisvert and C. W. Clark, (Eds.), NIST handbook of mathematical functions, Cambridge university press, 2010.

\bibitem{25}
W. B. Gearhart and H. S. Shultz, The Function $sin x/x$, College Mathematics Journal, 21 (1990) (2), 90-99.



\end{thebibliography}
\end{document}